\documentclass[10pt,a4paper]{article}
\usepackage[T1]{fontenc}
\usepackage{lmodern}
\usepackage[margin=25mm]{geometry}
\usepackage[authoryear,round]{natbib}
\usepackage{amsmath,amssymb,amsthm}
\usepackage[hidelinks]{hyperref}
\newtheorem{theorem}{Theorem}[section]
\newtheorem{lemma}[theorem]{Lemma}
\newtheorem{proposition}[theorem]{Proposition}
\theoremstyle{remark}
\newtheorem{remark}[theorem]{Remark}
\newcommand{\fGIG}{\mathrm{fGIG}}
\newcommand{\MGIGF}{\mathrm{MGIG}^{\mathrm F}}
\newcommand{\FP}{\mathrm{FP}}
\newcommand{\calT}{\mathcal T}
\newcommand{\dd}{\mathrm d}
\DeclareMathOperator{\Tr}{Tr}
\DeclareMathOperator{\tr}{tr}

\newcommand{\papertitle}{An extended Matsumoto--Yor property for free GIG distributions}
\hypersetup{pdftitle={An extended Matsumoto--Yor property for free GIG distributions},
  pdfauthor={Shuhei Nasu and Noriyoshi Sakuma}}
\newcommand{\paperaffiliation}{Department of Mathematics, Graduate School of Science,\\
The University of Osaka, Toyonaka, Osaka 560-0043, Japan}
\newcommand{\paperabstract}{%
We prove a two-parameter extension of the Matsumoto--Yor property for free generalized inverse Gaussian distributions. The transformation is the reciprocal form of a known matrix Yang--Baxter map. Our proof combines the independence-preserving property for matrix generalized inverse Gaussian distributions with asymptotic freeness. Since the relevant spectra stay in compact subsets of the positive half-line, continuous functional calculus and polynomial approximation suffice to handle all inverses. The same argument yields the free Matsumoto--Yor property itself for all positive rates of the free Poisson law, including those for which this law has an atom at zero. We also discuss the limit in which one of the two parameters tends to zero: 
the limiting laws are free Poisson only when the rate is at least one, so for smaller rates this limit does not recover the free Matsumoto--Yor property.}

\begin{document}
\title{\papertitle}
\author{Shuhei Nasu \and Noriyoshi Sakuma}
\date{November 17, 2026}
\maketitle
\begin{center}\small\paperaffiliation\end{center}
\begin{abstract}\paperabstract\end{abstract}

\section{Introduction}

The Matsumoto--Yor property transforms an independent pair of generalized
inverse Gaussian (GIG) and gamma variables into another independent pair
\citep{MatsumotoYor,LetacWesolowski2000}.
\citet{CroydonSasada} found a two-parameter extension involving two GIG
distributions, related to discrete KdV systems.
\citet{SasadaUozumi} placed the underlying transformation among
quadrirational Yang--Baxter maps.
Its matrix independence property was proved by
\citet{LetacWesolowski2024}; the matrix Yang--Baxter equation and a
characterization were established by \citet{LetacPiccioniWesolowski}.

We prove the corresponding freeness-preserving property for free GIG
distributions. It extends the free Matsumoto--Yor property of
\citet{Szpojankowski}. The proof uses the reciprocal form of the known
matrix map, which is a linear combination of three inverses.
The only analytic point is to control these inverses when passing to
the random matrix limit. F\'{e}ral's spectral bounds allow a direct
polynomial approximation argument. Finally, we prove the free Matsumoto--Yor property itself, which corresponds to
    $\beta=0$, for all rates $\lambda>0$; it was previously known for $\lambda>1$
    \citep{Szpojankowski}. We also show that the limit $\beta\downarrow0$ of our family
    recovers the marginal laws of this property only when $\lambda\ge1$.

\section{The freeness-preserving property}

We use a tracial $C^*$-probability space $(\mathcal A,\varphi)$ with a
faithful state. Thus the spectrum of a self-adjoint element equals the
support of its distribution. For positive invertible elements, inverses are defined by
continuous functional calculus. Write $I$ for the unit and $X\sim\mu$ for the distribution of $X$;
freeness and convergence in mixed moments have their usual meanings
\citep{NicaSpeicher}.

For $\kappa\in\mathbb R$ and $p,q>0$, the free GIG law
$\fGIG(\kappa,p,q)$ has density
\begin{equation}\label{eq:fgig}
 \frac{1}{2\pi}\sqrt{(x-a)(b-x)}
 \left(\frac p x+\frac{q}{\sqrt{ab}\,x^2}\right)
 \mathbf1_{[a,b]}(x),
\end{equation}
where $0<a<b$ are uniquely determined by
\begin{equation}\label{eq:endpoints}
 1-\kappa+p\sqrt{ab}-\frac{q(a+b)}{2ab}=0,
 \qquad
 1+\kappa+\frac q{\sqrt{ab}}-\frac{p(a+b)}2=0.
\end{equation}
This is $\mu(p,q,\kappa)$ in \citet{HasebeSzpojankowski}.
If $X\sim\fGIG(\kappa,p,q)$ and $c>0$, changes of variables give
\begin{equation}\label{eq:scaling}
 X^{-1}\sim\fGIG(-\kappa,q,p),\qquad
 cX\sim\fGIG(\kappa,p/c,cq).
\end{equation}

For $\alpha,\beta>0$ and positive invertible $X,Y$, put
\begin{equation}\label{eq:map}
 \calT_{\alpha,\beta}(X,Y)=(U,V),\qquad
 \begin{aligned}
 U&=\frac\beta\alpha Y^{-1}
    +\left(1-\frac\beta\alpha\right)(\alpha X+Y)^{-1},\\
 V&=X^{-1}+(\beta-\alpha)(\alpha X+Y)^{-1}.
 \end{aligned}
\end{equation}
These are positive invertible self-adjoint elements, even when a
coefficient in \eqref{eq:map} is negative. Indeed, with $S=\alpha X+Y$,
the order inequalities $S^{-1}\le Y^{-1}$ and
$\alpha S^{-1}\le X^{-1}$ imply
\begin{equation}\label{eq:positive}
 U\ge\min\{1,\beta/\alpha\}Y^{-1},\qquad
 V\ge\min\{1,\beta/\alpha\}X^{-1}.
\end{equation}

\begin{theorem}\label{thm:main}
Let $\lambda\in\mathbb R$ and $\alpha,\beta,\gamma_1,\gamma_2>0$.
If $X,Y$ are free self-adjoint elements with
\begin{equation}\label{eq:input}
 X\sim\fGIG(-\lambda,\alpha\gamma_1,\gamma_2),\qquad
 Y\sim\fGIG(\lambda,\gamma_1,\beta\gamma_2),
\end{equation}
then $(U,V)=\calT_{\alpha,\beta}(X,Y)$ are free with
\begin{equation}\label{eq:output}
 U\sim\fGIG(-\lambda,\alpha\gamma_2,\gamma_1),\qquad
 V\sim\fGIG(\lambda,\gamma_2,\beta\gamma_1).
\end{equation}
\end{theorem}

For $\alpha=\beta$, the map is simply $(Y^{-1},X^{-1})$.
For $\alpha\ne\beta$, the theorem gives a nontrivial extension of the
free Matsumoto--Yor property.

\begin{remark}\label{rem:origin}
The map \eqref{eq:map} is already given in
\citet[Remark~2.2]{LetacPiccioniWesolowski}. Specifically, if
\[
 \phi^{(\alpha,\beta)}(x,y)=
 \bigl(y(I+\alpha xy)^{-1}(I+\beta xy),\,
       x(I+\beta yx)^{-1}(I+\alpha yx)\bigr)
\]
and $\theta(x,y)=(x,y^{-1})$, then
$\calT_{\alpha,\beta}=\theta\circ\phi^{(\alpha,\beta)}\circ\theta$.
The Yang--Baxter equation in that work concerns $\phi$; we do not assert
it for this reciprocal form.
\end{remark}

\section{Proof by polynomial approximation}

We first record the elementary continuity argument used in the proof.
The states and algebras in the next lemma may vary with $n$.

\begin{lemma}\label{lem:inverse}
For each $n\ge1$, let $(\mathcal A_n,\varphi_n)$ be a unital
$C^*$-probability space, and let $(\mathcal A,\varphi)$ be another such space. 
Denote the units of $\mathcal A_n$ and $\mathcal A$ by $I_n$ and $I$, respectively.
Let $A_n,B_n\in\mathcal A_n$ and $A,B\in\mathcal A$ be self-adjoint elements. Assume the following.
\begin{enumerate}
    \item There are constants $0<c<C<\infty$ such that, for all sufficiently large $n$,
      \begin{equation}\label{as1}
        cI_n\le A_n\le CI_n,\qquad cI_n\le B_n\le CI_n,
      \end{equation}
      and moreover $cI\le A\le CI$ and $cI\le B\le CI$.
    \item The pair $(A_n,B_n)$ converges in mixed moments to $(A,B)$, that is,
      $\varphi_n\bigl(P(A_n,B_n)\bigr)\to\varphi\bigl(P(A,B)\bigr)$ as $n\to\infty$
      for every noncommutative polynomial $P\in\mathbb C\langle x,y\rangle$.
\end{enumerate}
Then, for every fixed $\alpha>0$, the tuple $(A_n,B_n,A_n^{-1},B_n^{-1},(\alpha A_n+B_n)^{-1})$ converges in mixed moments, with respect to $\varphi_n$ and
$\varphi$, to $\bigl(A,B,A^{-1},B^{-1},(\alpha A+B)^{-1}\bigr)$.

If the assumptions on $B_n$ and $B$ are replaced by
$0\le B_n\le CI_n,\qquad 0\le B\le CI$, while all other assumptions remain unchanged, then $\bigl(A_n,B_n,A_n^{-1},(\alpha A_n+B_n)^{-1}\bigr)$ converges in mixed moments to $\bigl(A,B,A^{-1},(\alpha A+B)^{-1}\bigr)$.
\end{lemma}

\begin{proof}
Fix $\alpha>0$ and put
\[
 S_n=\alpha A_n+B_n,\qquad S=\alpha A+B,\qquad
 m=c\min\{1,\alpha\},\qquad M=(\alpha+1)C.
\]
The assumption \eqref{as1} and self-adjointness imply that
\[
 \sigma(A_n),\ \sigma(B_n),\ \sigma(S_n),
 \sigma(A),\ \sigma(B),\ \sigma(S)\subset[m,M],
\]
where $\sigma(T)$ denotes the spectrum of $T$.
For $0<\varepsilon\le1$, we can choose a polynomial $r_\varepsilon$ with real coefficients such that
$$
 \sup_{t\in[m,M]}\bigl|r_\varepsilon(t)-t^{-1}\bigr|<\varepsilon.
$$
For each $T\in\{A_n,B_n,S_n,A,B,S\}$, continuous functional calculus gives
\[
 \|r_\varepsilon(T)-T^{-1}\|
 =\sup_{t\in\sigma(T)}\bigl|r_\varepsilon(t)-t^{-1}\bigr|
 <\varepsilon.
\]
In particular,
\[
 \|T^{-1}\|\le m^{-1},\qquad
 \|r_\varepsilon(T)\|\le m^{-1}+\varepsilon\le m^{-1}+1.
\]

Fix a mixed moment of length $k\ge1$, written as $\varphi_n(W_n)$, where
\[
 W_n=X_{n,1}\cdots X_{n,k},\qquad
 X_{n,j}\in\{A_n,B_n,A_n^{-1},B_n^{-1},S_n^{-1}\}.
\]
The choice of each factor is fixed as $n$ varies. Let $W=X_1\cdots X_k$
be the corresponding product formed from $A,B$.
Replace the inverse factors according to
\[
 A_n^{-1}\mapsto r_\varepsilon(A_n),\qquad
 B_n^{-1}\mapsto r_\varepsilon(B_n),\qquad
 S_n^{-1}\mapsto r_\varepsilon(S_n),
\]
leaving $A_n,B_n$ unchanged. Denote the resulting factors by
$\widetilde X_{n,j}$ and their product by
$\widetilde W_n=\widetilde X_{n,1}\cdots\widetilde X_{n,k}$.
Define $\widetilde X_j$ and $\widetilde W$ in the same way from $X_j$ and $W$.
With
\[
 L=\max\{C,m^{-1}+1\},
\]
we have, for every $n$ and $j$,
\[
 \|X_{n,j}-\widetilde X_{n,j}\|\le\varepsilon,
 \qquad
 \|X_{n,j}\|,\ \|\widetilde X_{n,j}\|\le L,
\]
and the same estimates hold for $X_j,\widetilde X_j$.

Replacing the factors one at a time gives the identity
\[
 W_n-\widetilde W_n
 =\sum_{j=1}^k
 \widetilde X_{n,1}\cdots\widetilde X_{n,j-1}
 (X_{n,j}-\widetilde X_{n,j})
 X_{n,j+1}\cdots X_{n,k},
\]
where an empty product means the unit. Hence
\[
 \|W_n-\widetilde W_n\|\le K\varepsilon,\qquad
 \|W-\widetilde W\|\le K\varepsilon,\qquad
 K=kL^{k-1}.
\]
In particular, $K$ is independent of $n$ and $0<\varepsilon\le1$.

For fixed $\varepsilon$, there is a noncommutative polynomial
$P_\varepsilon$ such that
\[
 \widetilde W_n=P_\varepsilon(A_n,B_n),\qquad
 \widetilde W=P_\varepsilon(A,B).
\]
Convergence in mixed moments therefore implies
\[
 \varphi_n(\widetilde W_n)\longrightarrow\varphi(\widetilde W).
\]
Since states have norm one,
\begin{align*}
 |\varphi_n(W_n)-\varphi(W)|
 &\le \|W_n-\widetilde W_n\|+\|\widetilde W-W\|\\
 &\quad+|\varphi_n(\widetilde W_n)-\varphi(\widetilde W)|\\
 &\le 2K\varepsilon
     +|\varphi_n(\widetilde W_n)-\varphi(\widetilde W)|.
\end{align*}
It follows that
\[
 \limsup_{n\to\infty}|\varphi_n(W_n)-\varphi(W)|
 \le 2K\varepsilon.
\]
Letting $\varepsilon\downarrow0$ proves convergence of every mixed moment.

Finally, suppose that
$0\le B_n\le CI_n$ and $0\le B\le CI$. Then
\[
 \alpha cI_n\le S_n\le MI_n,\qquad
 \alpha cI\le S\le MI,
\]
so $A_n,S_n,A,S$ still have spectra in $[m,M]$.
For a product $W_n$ with factors in
$\{A_n,B_n,A_n^{-1},S_n^{-1}\}$, replace only
$A_n^{-1}$ and $S_n^{-1}$ by $r_\varepsilon(A_n)$ and
$r_\varepsilon(S_n)$, respectively, and leave $A_n,B_n$ unchanged.
The estimates
\[
 \|B_n\|,\ \|B\|\le C\le L,\qquad
 \|r_\varepsilon(T)-T^{-1}\|<\varepsilon
 \quad(T\in\{A_n,S_n,A,S\})
\]
show that all factor bounds above remain valid.
The product identity and polynomial moment convergence thus give
\[
 \limsup_{n\to\infty}|\varphi_n(W_n)-\varphi(W)|
 \le 2kL^{k-1}\varepsilon
\]
for every such product. Letting $\varepsilon\downarrow0$ proves
the last assertion.
\end{proof}

Let $\Omega_N^+$ be the cone of real symmetric positive definite
$N\times N$ matrices and set $\tr_N=N^{-1}\Tr$.
Write $Z_N\sim\MGIGF_N(\kappa,p,q)$ when its density on this cone is
\begin{equation}\label{eq:matrix-density}
 f_N(z)\propto
 (\det z)^{N\kappa/2-(N+1)/2}
 \exp\left\{-\frac N2\Tr(pz+qz^{-1})\right\}.
\end{equation}
The factor $N/2$ is the normalization for real symmetric matrices.
The eigenvalue density has Vandermonde exponent $1$ and external field
\[
 V_N(x)=(1+N^{-1}-\kappa)\log x+px+q/x
 \quad\text{in the factor }\exp\{-\tfrac N2\textstyle\sum_i V_N(x_i)\}.
\]
Consequently, \citet[Theorems~1--3 and the subsequent real-case
remark]{Feral} gives, almost surely,
\begin{equation}\label{eq:spectral}
 \frac1N\sum_{j=1}^N\delta_{\lambda_j(Z_N)}
 \Longrightarrow\fGIG(\kappa,p,q),\qquad
 \lambda_{\min}(Z_N)\to a,
 \quad\lambda_{\max}(Z_N)\to b,
\end{equation}
with $a,b$ from \eqref{eq:endpoints}. These statements also hold when
$(\kappa,p,q)$ varies with $N$ and converges to a triple with $p,q>0$.
In particular, the spectra eventually lie in a fixed compact subset of
$(0,\infty)$.

Independent matrices of the form \eqref{eq:matrix-density} are
orthogonally invariant. By orthogonal Haar asymptotic freeness
\citep[Theorem~5.2]{CollinsSniady}, they therefore converge almost surely
in mixed moments to free elements with their respective free GIG laws.
To check the hypotheses, condition on the eigenvalues and use independent
Haar orthogonal eigenvector matrices. The increasingly ordered diagonal
eigenvalue matrices have mixed moment limits by bounded convergence of
their quantile functions, and are unchanged by transposition.
The same argument applies to independent orthogonally invariant matrices
with deterministic compactly supported marginal limits and eventual
uniform spectral bounds, including singular Wishart matrices.

Finally, the matrix independence property of
\citet[Theorem~5.2]{LetacWesolowski2024}, in the reciprocal form of
\citet[Theorem~1.1 and Remark~2.2]{LetacPiccioniWesolowski}, states that
independent inputs
\begin{equation}\label{eq:matrix-input}
 X_N\sim\MGIGF_N(-\lambda,\alpha\gamma_1,\gamma_2),\qquad
 Y_N\sim\MGIGF_N(\lambda,\gamma_1,\beta\gamma_2)
\end{equation}
give independent outputs $(U_N,V_N)=\calT_{\alpha,\beta}(X_N,Y_N)$ with
\begin{equation}\label{eq:matrix-output}
 U_N\sim\MGIGF_N(-\lambda,\alpha\gamma_2,\gamma_1),\qquad
 V_N\sim\MGIGF_N(\lambda,\gamma_2,\beta\gamma_1).
\end{equation}
In the notation of the latter reference, take the shape parameter
$-N\lambda/2$ and matrix parameters
$a=(N\gamma_1/2)I_N$, $b=(N\gamma_2/2)I_N$, and invert the second
input and output.

\begin{proof}[Proof of Theorem~\ref{thm:main}]
Let $X_N$ and $Y_N$ be independent random matrices
with the distributions in \eqref{eq:matrix-input}.
By the preceding asymptotic freeness argument, their mixed moments converge almost surely to those of $(X,Y)$, since the mixed moments of a free independent pair is determined by its marginal distributions.
Fix $0<c<C$ such that the open interval $(c,C)$ contains the supports of the two laws in \eqref{eq:input}. 
By \eqref{eq:spectral}, almost surely the spectra of
$X_N$ and $Y_N$ are contained in $[c,C]$ for all sufficiently large $N$; the spectra of $X$ and $Y$ are contained in $[c,C]$ as well.

Lemma~\ref{lem:inverse} thus gives, for every noncommutative polynomial $P$,
\[
 \tr_N P(U_N,V_N)\longrightarrow\varphi(P(U,V))
 \quad\text{almost surely}.
\]
On the other hand, \eqref{eq:matrix-output} and the same asymptotic
freeness result identify this limit as the mixed moments of two
freely independent self-adjoint elements with the respective
distributions in \eqref{eq:output}. Equality of the two limits
therefore shows that $U$ and $V$ have the respective distributions
in \eqref{eq:output} and are freely independent.
\end{proof}

\begin{remark}\label{rem:gap}
The positive lower spectral bound is needed for the matrix approximation,
not merely to define the limiting inverses. Weak convergence of empirical
laws to a measure supported away from zero does not exclude small
outlying eigenvalues. Here the endpoint convergence in \eqref{eq:spectral}
supplies precisely this bound. No joint operator-norm convergence is used.
Theorem~\ref{thm:main} also holds for a nonfaithful tracial state if $X,Y$
are explicitly assumed positive and invertible: apply the same polynomial
approximation on intervals containing their actual spectra.
\end{remark}

\section{The free Matsumoto--Yor property and limiting cases}

Let $\FP(\lambda,c)$ denote the free Poisson law of rate $\lambda>0$
and scale $c>0$:
\[
 \FP(\lambda,c)=(1-\lambda)_+\delta_0+
 \frac{\sqrt{(b-x)(x-a)}}{2\pi cx}\mathbf1_{[a,b]}(x)\,\dd x,
\]
where $a=c(1-\sqrt\lambda)^2$ and $b=c(1+\sqrt\lambda)^2$.
The usual free Matsumoto--Yor property was proved for $\lambda>1$ by
\citet[Theorem~3.3]{Szpojankowski}. The argument above also covers
the singular range, as follows.

\begin{proposition}\label{prop:my}
Let $\lambda,p,q>0$. If $X,Y$ are free self-adjoint elements with
$X\sim\fGIG(-\lambda,p,q)$ and $Y\sim\FP(\lambda,1/p)$, then
\[
 (X+Y)^{-1}\quad\text{and}\quad X^{-1}-(X+Y)^{-1}
\]
are free with laws $\fGIG(-\lambda,q,p)$ and $\FP(\lambda,1/q)$,
respectively.
\end{proposition}

\begin{proof}
Choose integers $k_N\ge1$ with $k_N/N\to\lambda$ and independent matrices
\[
 X_N\sim\MGIGF_N(-k_N/N,p,q),\qquad
 Y_N=(Np)^{-1}G_NG_N^{\mathsf T},
\]
where $G_N$ is an $N\times k_N$ standard real Gaussian matrix.
The matrix Matsumoto--Yor property \citep{LetacWesolowski2000}, valid
also for singular Wishart laws
\citep[Eq.~(1.3)]{LetacPiccioniWesolowski}, gives independent outputs
$U_N=(X_N+Y_N)^{-1}$ and $V_N=X_N^{-1}-(X_N+Y_N)^{-1}$ with
\[
 U_N\sim\MGIGF_N(-k_N/N,q,p),\qquad
 V_N\stackrel{\mathrm d}= (Nq)^{-1}G_NG_N^{\mathsf T}.
\]
Here the Wishart shape $k_N/2$ is admissible also when $k_N<N$;
the displayed equality specifies only the marginal law.
F\'{e}ral's theorem and the Wishart spectral limits
\citep{MarchenkoPastur,BaiYin} give the marginal limits and eventual
uniform upper bounds for both pairs. Moreover, $X_N\ge cI_N$ eventually
almost surely for some $c>0$. Thus $X_N$ and $X_N+Y_N$ have uniformly
bounded inverses, even if $Y_N$ is singular.
Apply orthogonal asymptotic freeness to the input and output pairs and
use the last assertion of Lemma~\ref{lem:inverse}. Comparing the two
mixed moment limits as above proves the result.
\end{proof}

It is useful to distinguish this property from the limit of the
two-GIG family. By \citet[Proposition~2.13]{HasebeSzpojankowski},
as $q\downarrow0$,
\begin{equation}\label{eq:boundary}
 \fGIG(\kappa,p,q)\Longrightarrow L_{\kappa,p}:=
 \begin{cases}
 \FP(\kappa,1/p),&\kappa\ge1,\\[2pt]
 \displaystyle\frac{1-\kappa}{2}\delta_0+
 \frac{1+\kappa}{2}\FP\left(1,\frac{1+\kappa}{2p}\right),
       &-1<\kappa<1,\\[4pt]
 \delta_0,&\kappa\le-1.
 \end{cases}
\end{equation}
In Theorem~\ref{thm:main}, this gives
$Y_\beta\Longrightarrow L_{\lambda,\gamma_1}$ and
$V_\beta\Longrightarrow L_{\lambda,\gamma_2}$ as $\beta\downarrow0$.
The inverse term in \eqref{eq:map} requires care: by \eqref{eq:scaling},
\begin{equation}\label{eq:inverse-boundary}
 \beta Y_\beta^{-1}\sim\fGIG(-\lambda,\gamma_2,\beta\gamma_1)
 \Longrightarrow L_{-\lambda,\gamma_2}.
\end{equation}
Thus a coefficient tending to zero does not by itself make this term
negligible. Its limiting law is $\delta_0$ precisely when $\lambda\ge1$;
this is a statement about distributions, not operator norms.

At $\beta=0$, define the map directly by
\begin{equation}\label{eq:beta-zero}
 U=(\alpha X+Y)^{-1},\qquad
 V=X^{-1}-\alpha(\alpha X+Y)^{-1}.
\end{equation}
No inverse of $Y$ is needed. If
$X\sim\fGIG(-\lambda,\alpha\gamma_1,\gamma_2)$ and
$Y\sim\FP(\lambda,1/\gamma_1)$ are free, Proposition~\ref{prop:my}
applied to $\alpha X,Y$ gives free outputs with laws
$\fGIG(-\lambda,\alpha\gamma_2,\gamma_1)$ and
$\FP(\lambda,1/\gamma_2)$ for every $\lambda>0$.
For $\lambda\ge1$, these laws agree with the marginal limits of the
two-GIG family. For $0<\lambda<1$, that family instead has the mixture
limits in \eqref{eq:boundary}. This difference is not a failure of the
free Matsumoto--Yor property.

The other limit, $\alpha\downarrow0$, is most simply described through reciprocals:
\[
X_\alpha^{-1}\sim\fGIG(\lambda,\gamma_2,\alpha\gamma_1)
 \Longrightarrow L_{\lambda,\gamma_2}.
\]
For $\lambda\ge1$, the limiting law has no atom at zero, so
$X_\alpha$ converges weakly to the reciprocal of
$\FP(\lambda,1/\gamma_2)$; at $\lambda=1$ this reciprocal law has
unbounded support. For $-1<\lambda<1$, mass $(1-\lambda)/2$ escapes
to infinity, and for $\lambda\le-1$ all mass does.
These are assertions about measures; bounded invertibility of a limiting
operator is not implied when its spectral support reaches zero.

\section*{Acknowledgements}

The authors thank Jacek Weso{\l}owski for helpful discussions
and comments on this work.

N. Sakuma was supported in part by JSPS KAKENHI Grant Numbers
JP25H00593, JP23K03133, and JP26K06828.
S. Nasu was supported by JST SPRING, Grant Number JPMJSP2138.

\section*{Use of generative AI}

During the preparation of this manuscript, the authors used
ChatGPT (OpenAI) to assist with English editing, improving the
clarity of mathematical explanations, and checking references
and manuscript consistency. The authors reviewed and revised
the AI-assisted material and take full responsibility for the
content of the article.

\end{document}